\documentclass[11pt,a4paper]{amsart}

\usepackage{pstricks,pst-node,pst-plot}
\usepackage{a4,amssymb,latexsym}

\title[Minimal Lagrangian 2-tori in $\CP^2$.]%
{Minimal Lagrangian 2-tori in $\CP^2$ come in real families of every dimension.}
\author{Emma Carberry}
\address{Department of Mathematics\\ Massachusetts Institute of Technology\\
77 Massachusetts Avenue \\ Cambridge, MA 02139-4307, USA}
\email{carberry@math.mit.edu}
\author{Ian McIntosh}
\address{Department of Mathematics \\ University of York \\ Heslington, York
YO10 5DD, UK}
\email{im7@york.ac.uk}
\thanks{This research was partially supported by an LMS Scheme 4 grant no.\ 
4720, and the first author was also partially supported by NSF grants DMS 9802479 and DMS 0306947.}

\newcommand{\C}{\mathbb{C}}
\newcommand{\R}{\mathbb{R}}
\renewcommand{\P}{\mathbb{P}}
\newcommand{\CP}{\mathbb{CP}}

\newcommand{\caX}{\mathcal{X}}
\newcommand{\eps}{\varepsilon}
\newcommand{\bft}{\mathbf{t}}

\renewcommand{\o}{\omega}
\newcommand{\om}{\omega^m}
\newcommand{\x}{\ensuremath{\xi_{1}}}
\newcommand{\xx}{\ensuremath{\xi_{2}}}
\newcommand{\y}{\eta_{1}}
\newcommand{\yy}{\eta_{2}}
\newcommand{\xm}{\ensuremath{\xi^m_{1}}}
\newcommand{\xxm}{\ensuremath{\xi^m_{2}}}
\newcommand{\ym}{\eta^m_{1}}
\newcommand{\yym}{\eta^m_{2}}
\newcommand{\thj}{\theta_j}
\newcommand{\thmj}{\theta_{mj}}
\newcommand{\2}{\frac{2\pi}{3}}
\newcommand{\xmj}{X_{mj}}
\newcommand{\ymj}{Y_{mj}}
\newcommand{\z}{\zeta}
\newcommand{\cm}{c_m}
\newcommand{\f}{\frac}
\newcommand{\res}{\mbox{Res}}

\newtheorem{theor}{Theorem}
\newtheorem{prop}{Proposition}
\newtheorem{lem}{Lemma}

\begin{document}
\begin{abstract}
We show that for every non-negative integer $n$, there is a real
$n$-dimensional family of minimal Lagrangian tori in $\CP^2$, and hence of
special Lagrangian cones in $\C^3$ whose link is a torus. The proof utilises
the fact that such tori arise from integrable systems, and can be described
using algebro-geometric (spectral curve) data.
\end{abstract}

\maketitle

\section{Introduction.}
Minimal Lagrangian tori in $\CP^2$ are of some interest at present, since they
classify the special Lagrangian cones in $\C^3$ whose link is a 2-torus. 
It is quite straightforward to show that a minimal Lagrangian surface in $\CP^2$
lifts locally to a minimal Lagrangian cone in $\C^3$, which must be congruent to a
special Lagrangian cone by a result of Harvey and Lawson \cite{HL}. 
This lift essentially extends globally: if the surface itself 
does not lift then a three-fold cover of it does (for tori cf.\ \cite{BPW}). 
Conversely, every special 
Lagrangian (regular) cone in $\C^3\setminus\{0\}$ projects to a minimal 
Lagrangian surface in $\CP^2$. 
Now, one can argue (as Joyce does in \cite{Joy02}) that
special Lagrangian cones in $\C^3$ are the local models for the 
singularities which may occur in special Lagrangian fibrations of Calabi-Yau 3-folds. 
Such fibrations form the core of the Strominger-Yau-Zaslow conjecture concerning the source of
mirror symmetry \cite{StrYZ}. Since the smooth fibres in these fibrations are 3-tori one
expects that understanding cones over tori will give some insight into the singular
fibres. Many examples have already been constructed (see, for example, 
\cite{Joy01}), but
it would be useful to know how rich the collection of all such cones is. 
Our results will show that the
moduli space of minimal Lagrangian tori in $\CP^2$ has a very complicated
structure, including countably many connected components of every real dimension. 
To do this we need to recall what is already known about minimal tori in
$\CP^2$.

Minimal surfaces in $\CP^2$ come in two types, known as complex isotropic (or
super-minimal) and non-isotropic. All minimal Lagrangian tori are 
non-isotropic.  For some time now it has been known 
that the construction of non-isotropic minimal tori 
reduces to a periodicity problem in the
theory of algebraically integrable Hamiltonian systems. 
There is a bijective correspondence between non-isotropic minimal tori in 
$\CP^2$ and spectral data $(X,\lambda,\mathcal{L})$ where $X$ is an algebraic 
curve, $\lambda:X\to\CP^1$ is a three-sheeted branched covering,  and $\mathcal{L}$ 
is a line bundle on $X$, all of which satisfy certain conditions (see \cite{McI03} 
for details). For ease of discussion, we will call $(X,\lambda)$ the spectral curve
and the genus of $X$ the spectral genus. 
For minimal Lagrangian tori we know from \cite{McI03} that the spectral genus is even 
and the spectral curve possesses an additional holomorphic involution. 
Our aim is to show that
there exist minimal tori for every possible spectral genus, and minimal Lagrangian tori
for every even spectral genus.

The approach we will follow, which is by now standard, is to consider 
the tori as lying in the much larger set of all non-isotropic minimal immersions of 
$\R^2$ which have the property of being ``of finite type''. 
One knows from \cite{McI95} that for every
spectral genus there are minimal immersions of $\R^2$ 
into $\CP^2$ of finite type for a smooth family of spectral curves,
but the issue of double periodicity of these maps is a delicate question 
concerning the rationality of a certain real 
2-plane with respect to the period lattice of a (generalised) Jacobi variety
determined by the spectral curve.
Our aim is to show that there are, for each spectral genus,
spectral curves for which this 2-plane is rational.
\begin{theor} 
\label{th:every_genus}
For every integer $g\geq 1$ there are countably many spectral curves 
of genus $g$ that give rise to (non-congruent) minimal immersions of tori into 
$\CP^2$. Further, when $g$ is even, there are countably many spectral curves of 
genus $g$ that give rise to (non-congruent) minimal Lagrangian immersions of tori into $\CP^2$.
\end{theor} 
For $g=0$ there is precisely one possible choice of spectral
curve (namely, the Riemann sphere with the three-fold cover 
$\zeta\mapsto\zeta^3$).
This gives one well-known minimal Lagrangian torus, which is an $S^1\times S^1$
orbit: this is the most symmetric of possibilities.
For minimal Lagrangian tori, the result of the theorem can be deduced for $g=2$ and $g=4$ 
from, respectively, Castro \& Urbano \cite{CasU} and
Joyce \cite{Joy01}. 

Whenever a spectral curve for a minimal immersion is given, it follows from 
\cite{McI01,McI03} that the line bundle $\mathcal{L}$ may be chosen from a 
real $g$~dimensional family, or in the Lagrangian case, from a real $g/2$~dimensional 
family. This choice also marks a point on the torus, so there is a two parameter family
of line bundles which yield the same minimal torus in $\CP^2$. 
Allowing for this we deduce, as a corollary, the second theorem.
\begin{theor}
\label{th:every_dim}
For each $g>2$, there are countably many real \mbox{$g-2$}~dimensional families 
of non-congruent minimal immersions of tori into $\CP^2$. Furthermore, whenever 
$g>4$ is even, there are countably many real \mbox{$g/2 - 2$}~dimensional families of 
non-congruent minimal Lagrangian immersed tori in $\CP^2$. Each family corresponds to 
a particular genus $g$ spectral curve. 
\end{theor}
These two theorems are not unexpected, following the discussion 
of this in \cite{McI96,McI03}. 
An analogous problem has been solved for constant mean curvature tori in $\R^3$ 
\cite{ErcKT,Jag} and for minimal tori in $S^3$ \cite{Car}. In both of those
cases the spectral curve is hyperelliptic and the problem was solved by 
studying the variation of
holomorphic and meromorphic differentials, which have a well-known formulation
in that case. 
For $\CP^2$ the spectral curves are trigonal (i.e., three-fold covers of the Riemann
sphere),
for which we were not aware of any directly relevant variational formulas, so 
the proof we have found requires some perserverance during the calculation 
of the variations of differentials. Fortunately we were able to
simplify this somewhat so that only variations of differentials on elliptic
curves needed to be computed.
The method we use should apply equally well to the same problem for 
minimal tori in
$\CP^n$, but this will not be a trivial calculation and we happily forsake it
for the specific case at hand. 

To summarise the method (and the paper) we first recall from 
\cite{McI96,McI03} that the problem reduces to the following. To each spectral
curve $(X,\lambda)$ we assign,
by fixing a basis for a certain space of meromorphic differentials, a real
2-plane $W_X$ in $\R^{g+2}$. By fixing the basis appropriately, the
periodicity problem is solved whenever $(X,\lambda)$ gives rise to a rational
2-plane $W_X\in Gr(2,g+2)$, the Grassmannian of real 2-planes in $\R^{g+2}$.

The spectral curve $X$ is essentially determined
by the branch divisor of $\lambda$. This is invariant under
the real involution $\lambda\mapsto\bar\lambda^{-1}$, and contains
the points $\lambda=0,\infty$ as branch points of index three. The
construction of minimal immersions requires that none of the branch points
lie on the circle $|\lambda|=1$. Nevertheless, the assignment of $W_X$ to $X$
still makes sense when a pair of branch points, $\lambda_j$ and 
$\bar\lambda_j^{-1}$ are pushed together onto the
unit circle. In that case $X$ acquires a node. In particular, we may consider the case where $X$
is rational with $g$ nodes on the unit circle. We construct a real $2g$
parameter family $\caX$ containing such a rational curve, in which every open
neighbourhood of this curve contains, generically, smooth spectral curves
(which necessarily give rise to minimal 2-planes). This gives us 
a smooth (in fact, real analytic) map
\[
W:\caX\to Gr(2,g+2)\ ;\ X\mapsto W_X
\]
between real $2g$~dimensional manifolds. Our aim is to show that $dW$ is
invertible at one of the rational curves and, by invoking the inverse function
theorem, deduce that $W$ is locally invertible. It follows that in some open
neighbourhood of this curve there are countably many
spectral curves for which $W_X$ is rational.

We construct the family $\caX$ and compute $dW$ by constructing, for each
node on the unit circle, a partial desingularisation through a one parameter
family of nodal elliptic curves; in effect we ``pull apart'' that node to get a
pair of branch points. 

To restrict our attention to minimal Lagrangian tori we must consider spectral
curves of even
genus $g=2p$ whose branch divisor is also invariant under the
holomorphic involution $\lambda\mapsto -\lambda$. This condition specifies a
real $2p$~dimensional subfamily $\caX_L\subset\caX$ for which $W_X$ lies in a
real \mbox{$p+2$}~dimensional subspace of $\R^{g+2}$; by a judicious choice of
coordinates we identify this subspace with $\R^{p+2}$. We therefore obtain,
by restriction, a smooth map
\[
W:\caX_L\to Gr(2,p+2).
\]
By a simple adaptation of the argument for minimal tori, we show that $dW$ is
invertible at one nodal curve and again deduce that nearby there are countably many
spectral curves for minimal Lagrangian tori.

\section{Families of spectral curves.}

\subsection{Spectral curves and the periodicity conditions.}
A compact Riemann surface $X$ of genus $g$ is the spectral curve for a minimal
immersion of $\R^2$ into $\CP^2$ if it 
admits a degree three cover $\lambda:X\to\C_\infty$ of the Riemann sphere
whose branch divisor has
the form
\[
B = 2.0+\lambda_1+\ldots+\lambda_{2g}+2.\infty
\]
which is invariant under the real involution $\lambda\mapsto \bar\lambda^{-1}$
and where none of the branch points $\lambda_j$ lie on the unit circle. It is not
necessary that the $\lambda_j$ be distinct but we will only work
with this case here. Given such a branch divisor $B$ one can construct a
spectral curve $X$ by gluing together three copies of $\C_\infty$ along
branch cuts which connect $0$ with $\infty$ and which on two sheets have
branch cuts between $\lambda_j$ and $\bar \lambda_j^{-1}$ (see figure 1). 
Then $X$ comes equipped with a unique real involution $\rho$ for which 
$\rho^*\lambda
= \bar\lambda^{-1}$ and $\rho$ fixes every point lying over $|\lambda|=1$.
\begin{figure}
\begin{center}
\begin{pspicture}(-3,-3)(3,3)
\pscircle[linestyle=dashed](0,0){2}
\psline{*-}(0,0)(3,0)
\psline{*-*}(1,1)(2,2)
\uput[0](1,1){$\lambda_j$}
\uput[0](2,2){$\bar\lambda^{-1}_j$}
\psline{*-*}(-1,0)(-2.8,0)
\rput(-1,-0.3){$\lambda_k$}
\rput(-2.8,-0.3){$\bar\lambda^{-1}_k$}
\end{pspicture}
\end{center}
\caption{}
\end{figure}
Let $O_1,O_2,O_3$ be the distinct points on $X$ lying over $\lambda=1$,
and let $H^0(\Omega^\prime_X)$ denote the space of meromorphic
differentials on $X$ with no poles except possibly simple ones at these
points. This space is complex \mbox{$g+2$}~dimensional with 
real subspace
\[
V_X=\{\omega\in H^0(\Omega^\prime_X):\rho^*\omega = -\bar\omega\}.
\]
The dual space $V_X^*$ has a natural lattice $\Lambda$ consisting of 
integrals over
homology cycles $A$ for which $\rho(A)$ is homologous to 
$A^{-1}$.
We can see (from figure 1) that $\Lambda$ is generated by cycles
$A_1,\ldots,A_{g+2}$ represented in the cut plane
by (for $j=1,\ldots,g$) 
cycles encircling only a pair of branch points $\lambda_j,\bar \lambda_j^{-1}$ or (for 
$j=g+1,g+2$) a cycle encircling $O_1$ or $O_2$. These generators form a
basis for $V_X^*$ with which we identify it with $\R^{g+2}$. 

Now let $P_0$ be the unique point on $X$ lying over
$\lambda=0$. Fix a local parameter $\zeta$ about $P_0$ such that
$\zeta^3=\lambda$ and define, for each $\omega\in H^0(\Omega^\prime_X)$,
\[
\omega(P_0) = Res_{P_0}\zeta^{-1}\omega.
\]
This provides two real linear maps $u,v\in V_X^*$ given by taking the real and
imaginary parts of $\omega(P_0)$:
\[
\omega(P_0)= u(\omega) +iv(\omega). 
\]
The geometry behind these definitions is as follows (cf.\ \cite{McI01}). The
space $H^0(\Omega^\prime_X)$ is the space of regular differentials on the
singular curve $X^\prime$ obtained from $X$ by identifying the three points
$O_1,O_2,O_3$ together to obtain a triple node $O$. 
The quotient of its dual
space by the homology lattice of $X^\prime\setminus\{O\}$ gives the generalised
Jacobian $J(X^\prime)$ of $X^\prime$. Inside this non-compact complex Lie
group lies the real compact subgroup $J_R(X^\prime)\simeq V_X^*/\Lambda$ of
fixed points for the real automorphism $\rho^*$ covered by $\omega\mapsto
-\overline{\rho^*\omega}$. 

In \cite{McI96, McI01} it is proven that any minimal $f:\R^2\to \CP^2$ 
of finite type factors through $J_R(X^\prime)$, 
i.e. $f=\Theta\circ  L^\prime$, where
\[
\R^2\stackrel{L^\prime}\to J_{R}(X^\prime) \stackrel{\Theta}\to \CP^2,
\]
and moreover that $f$ is doubly-periodic precisely when $L^\prime$ is. 
The map $L^\prime$ is linear, and hence is determined by its differential 
at $0$. This (see  \cite{McI01}) satisfies
\[
dL^\prime_0\left(\frac{\partial}{\partial z}\right) = 
d\mathcal{A}^\prime_{P_0}\left(\frac{\partial}{\partial\zeta}\right),
\]
where $z$ is the usual coordinate on $\C$ and $\mathcal{A}^\prime_{P_0}$ 
is the Abel-Jacobi map of $X^\prime$ with basepoint $P_0$. It follows that
$dL^\prime_0\left(\frac{\partial}{\partial z}\right)$ is given by
\[
\omega\mapsto\frac{\partial}{\partial\zeta}\int_{P_0}^\zeta \omega = \omega(P_0)
\]
so that 
\[
dL^\prime(T_0\R^2) = Span_\R (u, v) \subset V_X^*.
\]
Let $W_X\subset \R^{g+2}$ denote this real 2-plane. 
\begin{theor}[\cite{McI96}]
A minimal immersion of $\R^2$ into $\CP^2$ with spectral data $(X,\lambda)$
is doubly periodic if and only if $W_X\in Gr(2,g+2)$ is a rational 2-plane.
\end{theor}
Notice that all that is required to obtain $W_X$ is that $X$ be a three-fold
cover of the Riemann sphere, satisfying the reality condition, which has a
branch of order three over $\lambda=0$ and 
no branch over $\lambda=1$. In particular, $W_X$ is two dimensional even
when the branch points degenerate in pairs to nodes over the unit circle away
from $\lambda=1$. These are the degenerations we intend to use.

\subsection{A smooth family of spectral curves.}
Because there
are finitely many spectral curves with the same branch locus, we must be
more precise about the spectral data $(X,\lambda)$ which we assign to each
branch divisor and do it in a way which aids the calculation
in mind. We only care about constructing some family of spectral curves
which degenerate to a rational $g$-nodal curve.

To this end, let ${\P}_\zeta$ denote the Riemann sphere equipped with rational
coordinate $\zeta$, and think of this as a 3-fold cover of $\C_\infty$ by the
map $\zeta\to\zeta^3$. We start by fixing $g$ distinct points $b_1,\ldots,b_g$ 
on the unit circle in $\C_\infty$, and choose cube roots $c_1,\ldots,c_g$,
$c_m^3=b_m$, which we also write in the form $c_m=e^{i\theta_m}$, 
and think of as lying on ${\P}_\zeta$. We define
$X_b$ to be the rational $g$-nodal curve obtained from ${\P}_\zeta$ by 
identifying each pair $c_m,\eps c_m$, where $\eps = \exp(2\pi i/3)$. Now we
construct, for each $j$, a 1-parameter family $X_b(t_j)$ of 
elliptic curves with $g-1$ nodes for which $X_b(0)=X_b$.

For each $j$ we define the real 1-parameter family of elliptic 
curves $E_j(t_j)$ given by the equation
\begin{equation} 
\frac{\lambda}{b_j}+\frac{b_j}{\lambda} = 2p_{t_j}(x), \label{eq:lambda}
\end{equation}
where
\[
p_{t_j}(x)=(x-2a)^2(x-t_j+a)-1,\text{ for } a=2^{-1/3}, 
\]
and where $t_j$ is a real parameter. For $t_j\neq 0$
this curve is non-singular and the map $\lambda:E_j\to\C_\infty$;
$(x,\lambda)\mapsto\lambda$ exhibits $E_j$ as a spectral curve of genus 1. 
Note that
\begin{eqnarray*}
p_{t_j}(x)^2-1 &=& \frac{1}{4}\left(\f{\lambda}{b_j}-\f{b_j}{\lambda}\right)^2\\
		&=& \left(\f{\lambda}{b_j}-p_{t_j}(x)\right)^2,
\end{eqnarray*}
so
\begin{equation*}
(x-t_j+a)(p_{t_j}(x) -1) = \left(\frac{\lambda b_j^{-1}-p_{t_j}(x)}{x-2a}\right)^2.
\end{equation*}
Setting
\begin{equation*}
y=\frac{\lambda b_j^{-1}-p_{t_j}(x)}{x-2a} 
\end{equation*}
thus realises $E_j(t_j)$ as the immersed curve 
\begin{equation}
\label{eq:y}
y^2=(x-t_j+a)(p_{t_j}(x)-1)
\end{equation}
in $\CP^2$.

For
$t_j=0$, $E_j$ becomes a rational nodal curve with rational parameterisation
\begin{equation}
{\P}_\zeta\to E_j(0);\ x(\zeta)=\frac{a\eps}{c_j\zeta}(\zeta-c_j)(\zeta-\eps
c_j),\ \lambda(\zeta)=\zeta^3. \label{eq:zeta}
\end{equation}

It follows from this that $E_j(0)$ is the nodal curve obtained from ${\P}_\zeta$
by identifying $c_j$ with $\eps c_j$.
On $E_j(t_j)$ there are, for each $m=1,\ldots,
 g$, two points $P_m,Q_m$ for which
$\lambda(P_m)=\lambda(Q_m)=b_m$, distinguished by the property that at $t_j=0$
these have $\zeta(P_m)=c_m$ and $\zeta(Q_m)=\eps c_m$. We define $X_b(t_j)$ 
to be the singularisation of $E_j(t_j)$ obtained by identifying each pair
$P_m,Q_m$ for $m\neq j$. Clearly $X_b(0)=X_b$. The points $O_k(t_j)$ lying
over $\lambda=1$ are labelled 
so that at $t_j=0$ they satisfy $\zeta(O_1)=1$, $\zeta(O_2)=\eps$ and $\zeta(O_3)=\eps^2$.

This collection of curves $X_b(t_j)$ uniquely determines, for each
$\bft=(t_1,\ldots,t_g)$ (assuming each $t_j$ sufficiently small and non-zero) a
non-singular genus $g$ spectral curve $X_b(\bft)$ with the property that its
branch divisor has, for each $j=1,\ldots g$, the same pair $\lambda_j,\bar
\lambda_j^{-1}$ as $X_b(t_j)$ with the same sheet identifications in the
construction as a 3-sheeted cover of $\C_\infty$. Now we consider the same
construction for each point in a family
of deformations of the node locations, varying the angular parameters
$\theta_j$ while keeping $0<\theta_1<\ldots<\theta_g<2\pi/3$.
In this way we obtain a real $2g$-parameter family $\caX_b$ of, 
at worst nodal, Riemann surfaces $X_b(\bft,\theta)$ parameterised by 
\[
\bft=(t_1,\ldots,t_g)\ ,\theta=(\theta_1,\ldots,\theta_g)
\]
with $X_b(0,0)=X_b$. To avoid complications we choose the domain of the
parameters so that
$\caX_b$ is diffeomorphic to an open $2g$-cube. By construction,
in every open neighbourhood of the origin there lies a spectral curve.

Over each Riemann surface $X$ in this family we have a basis
$\omega^1,\ldots,\omega^{g+2}$ for the real space $V_X$ 
of meromorphic differentials described earlier, normalised with respect to the
real homology cycles $A_1,\ldots,A_{g+2}$ so that 
\[
\oint_{A_j}\omega^m = \delta^m_j.
\]
We will, for simplicity, identify the linear maps $u,v$ on $V_X$ with their
coordinate representations in this basis, so that
\[
\omega^m(P_0) = u_m+iv_m,\ m=1,\ldots,g+2.
\]
We will denote by $u(\bft,\theta),v(\bft,\theta)$ the $\R^{g+2}$-valued functions
with these components $u_m,v_m$.
\begin{prop}
The map
\[
W:\caX_b\to Gr(2,g+2);\ W(\bft,\theta) =
Sp\{u(\bft,\theta),v(\bft,\theta)\}
\]
obtained from this construction is real analytic.
\end{prop}
\begin{proof} It suffices to prove that the complex values $\omega^m(P_0)$
vary analytically over this open cube $\caX_b$. 
Suppose $X\in\caX_b$. We can dissect $X$ into a union of the form
\[
X=S\cup N_1\cup\ldots\cup N_g
\]
where $S$ is a Riemann sphere with $2g$ discs excised and each $N_j$ is a
neck, possibly singular. Clearly, for $\caX_b$ suitably small,
we can choose $S$ to be the same for each
$X\in\caX_b$. Then $\zeta$ is a rational parameter on $S$ with
$\lambda=\zeta^3$. For each $(\bft,\theta)$, $\omega^m_{(\bft,\theta)}$ is
holomorphic on $S$. We write this as $f_m(\zeta,\bft,\theta)d\zeta$, so that
$\omega^m_{(\bft,\theta)}(P_0)=f_m(0,\bft,\theta)$. 
Since $\zeta=0$ is away from the moving branch points
we may deduce from the variational formulas of 
Fay \cite[Prop 3.7]{Fay} that $f_m(0,\bft,\theta)$ is a real analytic function.
\end{proof}

\section{Variation of differentials.}
We shall show that the differential of 
\[
W:\caX\to Gr(2,g+2)\ ;\ X\mapsto W_{X}
\]
is invertible at some $X_b$. If $e_1,e_2,\ldots,e_{g+2}$ 
is any basis for $\R^{g+2}$ with $e_{g+1},e_{g+2}$ a basis for $W_X$, then
\[
dW_X \mbox{ is invertible} \iff 
M=\left(\begin{array}{cc} a^i_j&b^i_j\\
c^i_j&d^i_j\end{array}\right) \mbox{ is invertible,}
\]
where $1\leq i,j\leq g$ and
\begin{eqnarray*}
\frac{\partial e_{g+1}}{\partial \theta_j}&=&e_i a^i_j + e_{g+1} a^{g+1}_j + e_{g+2} a^{g+2}_j,\\
\frac{\partial e_{g+2}}{\partial \theta_j}&=& e_i b^i_j +e_{g+1} b^{g+1}_j +e_{g+2} b^{g+2}_j,\\
\frac{\partial e_{g+1}}{\partial t_j}&=&e_i c^i_j + e_{g+1} c^{g+1}_j + e_{g+2} c^{g+2}_j,\\
\frac{\partial e_{g+2}}{\partial t_j}&=& e_i d^i_j +e_{g+1} d^{g+1}_j +e_{g+2} d^{g+2}_j .
\end{eqnarray*}
For each $m=1,\ldots,g+2$, pulling back ${\o}^m$ by the natural map  
\mbox{$E_{j}(t_{j})\to X_{b}(t_{j})$} gives a differential on $E_{j}(t_{j})$, 
which we also denote by ${\o}^{m}$, with the understanding that our calculations 
take place on $E_{j}(t_{j})$. For $m\neq j$, ${\o}^m$ is holomorphic except for 
simple poles at $P_{m}$ and $Q_{m}$ 
and is characterised by:
\begin{enumerate}
\item[(i)] ${\rm \res}_{P_{m} } {\o}^{m}=\frac{1}{2\pi i}$,
${\rm
\res}_{Q_{m} } {\o}^m=-\frac{1}{2\pi i}$, 
\item[(ii)] $\oint_{A_j}{\o}^m=0$, where $A_j$ is the real
homology cycle of $E_{j}(t_{j})$. 
\end{enumerate}
Here we take $P_{g+1}=O_1,\,Q_{g+1}=O_2,\,P_{g+2}=O_2$ and $Q_{g+2}=O_3$.
Clearly ${\o}^j$ is the holomorphic differential on $E_j(t_j)$
satisfying the normalisation \mbox{$\int_{A_j}{\o}^j = 1$.}

For $t_j = 0$, 
\[
\om=\frac{1}{2\pi i}\left(\frac{d\z}{\z-\cm} - 
\frac{d\z}{\z-\eps\cm}\right),
\]
where $c_{g+1}=1$ and  $c_{g+2}=\eps$. Then
\[
\om(P_0)=\frac{1}{2\pi i}\left(\f{1}{\eps\cm}-\f{1}{\cm}\right)=
\frac{\sqrt{3}}{2\pi}\exp[i(\frac{2\pi}{3}-\theta_m)],
\]
so
\begin{eqnarray}
u&=&\f{\sqrt{3}}{2\pi} \left(\cos(\2-\theta_1)
,\ldots, \cos(\2-\theta_g), -\f{1}{2}, 1 \right)\label{eq:u}\\
v&=&\f{\sqrt{3}}{2\pi} \left(\sin(\2-\theta_1)
,\ldots, \sin(\2-\theta_g), \f12, 0 \right)\label{eq:v}
\end{eqnarray}

In particular, we may take $e_m$  to be the standard basis vectors 
$e_m$ of $\R^{g+2}$ for $m=1,\ldots,g$ and $e_{g+1}=u$, $e_{g+2}=v$. 
Much of this section is devoted to the computation of the resulting matrix 
$M$; throughout $m$ shall be in the range $1,\ldots, g$.

\begin{prop} \label{thm:derivs}
\[
\begin{array}{rcl}
\left.{\partial u_m}/{\partial 
\theta_j}\right|_{\bft=0}&=&\left\{\begin{array}{ll}
\frac{\sqrt{3}}{2\pi}\sin(\2-\thj) &\mbox{ for } m= j\\
0&\mbox{ for }m\neq j \end{array}\right.\\
&\\
\left.{\partial v_m}/{\partial\theta_j}\right|_{\bft=0}&=
&\left\{\begin{array}{ll}
\frac{-\sqrt{3}}{2\pi}\cos(\2-\thj)&\mbox{ for } m= j\\
0&\mbox{ for }m\neq j \end{array}\right.\\
&\\
\left.{\partial u_m}/{\partial t_j}\right|_{\bft=0}&=&
\left\{\begin{array}{ll}
\frac{-\sqrt{3}}{2\pi a}\cos(\2-\thj)&\mbox{for }
m=j\\ 
&\\
\frac{1}{12\pi
a}\left[\cos(\2-\thj)\ymj+\sin(\2-\thj)\xmj
\right]&\mbox{for } 
m\neq j\end{array}\right.\\
&\\
\left.{\partial v_m}/{\partial t_j}\right|_{\bft=0}&=&
\left\{\begin{array}{ll}
\frac{-\sqrt{3}}{2\pi a}\sin(\2-\thj)&\mbox{for }
m=j\\ 
&\\
\frac{1}{12\pi
a}\left[\sin(\2-\thj)\ymj-\cos(\2-\thj)\xmj
\right]&\mbox{for } 
m\neq j\end{array}\right.
\end{array}                                
\]
where we have defined
\[
(\xm,\ym)=(x(P_m),y(P_m)),\; (\xxm,\yym)=(x(Q_m),y(Q_m))
\]
\[
\xmj=\frac{2a}{\xm}-\frac{2a}{\xxm},
\]
\[
\ymj =
\sin{3\theta_{mj}}\left(\frac{(\xm)^3-(\xxm)^3}{(\xm\xxm)^3}\right)
-\frac{\sqrt{3}}{2\cos{\thmj}+1},
\]
and
\[
\thmj=\theta_m - \thj.
\]
\end{prop}
The $\thj$ derivatives follow immediately from (\ref{eq:u}) and (\ref{eq:v}).
To obtain the $t_j$ derivatives we break the calculation into a number 
of lemmas. The first gives the case $m=j$. 
\begin{lem} 
\label{lem:omegajdot}
\[
\dot{{\o}^j}(P_0)=-\frac{\sqrt{3}}{2\pi a}\eps e^{-i\theta_j}
\]
where we have employed the notation 
$\dot{{\o}^m}(P_0)=\left.\frac{\partial}{\partial t_j}\right|_{\bft=0}{\o}^m(P_0)$ and so forth.
\end{lem}
\begin{proof}
Clearly $ {\o}^j=kdx/y$ for some $k=k(t_j)$.
On $E_j(t_j)$ set $\tilde x=1/x$ and $\tilde y=y/x^2$, 
then $(\tilde x(P_0),\tilde y(P_0))=(0,-1)$. Using (\ref{eq:lambda}), a calculation gives
\[
\lambda=\frac12 c_j^3\tilde x^3+O(\tilde x^4).
\]
From (\ref{eq:zeta}), at $t_j=0$,
\begin{eqnarray*}
\tilde x&=&\frac{c_j\zeta}{a\eps}(\zeta-c_j)^{-1}(\zeta-\eps c_j)^{-1}\\
&=&\frac{\eps}{ac_j}\zeta(1+O(\zeta)),
\end{eqnarray*}
so
\begin{equation}\label{eq:tilde_x}
\tilde x=\frac{\eps}{ac_j}(\zeta+O(\zeta^2))\mbox{ for all } t_j.
\end{equation}
This gives
\[
{\o}^j(P_0)=
\frac{\eps k}{ac_j}.
\]
For all $t_j$, ${\o}^j$ satisfies the normalisation $\oint_{A_j} {\o}^j =1$. 
Differentiating this and evaluating at $t_j=0$ gives
\[
k(0) \res_{P_j}\dot{\left(\frac{dx}{y}\right)} + 
\dot{k} \res_{P_j}\frac{dx}{y}=0.
\]
When $t_j=0$, $x(P_j)=0$ whilst $y^2=x^2(x-3a)(x+a)$, so 
\[
\res_{P_j}\frac{dx}{y}=\frac{1}{\sqrt{3}ai}
\]
and we obtain
\[
\frac{\sqrt{3}a}{2\pi} \res_{P_j}\dot{\left(\frac{dx}{y}\right)} + \frac{\dot{k}}{\sqrt{3}ai}=0.
\]
An easy calculation gives
\[
\dot{y}=\frac{-(x^3-3ax^2+1)}{y},
\]
and hence $\res_{P_j}\dot{\left(\frac{dx}{y}\right)}=\frac{2a}{\sqrt{3}i}$, so
\[
\dot{k}=\frac{-\sqrt{3}}{2\pi}.
\]
We conclude that
\[
\dot{{\o}^j}(P_0)=\frac{-\sqrt{3}\eps}{2\pi ac_j}.
\]
\end{proof}
Taking real and imaginary parts produces
the expressions in Proposition~\ref{thm:derivs}. 

We turn now to the case $m\neq j$. We shall fix such an $m$, and so dispense with various 
superscripts: we denote $\xi^m_k$ by $\xi_k$ and $\eta^m_k$ by $\eta_k$. 
The differential $\om$ has the form
\[
{\o}^m=\frac{\alpha(t_j) x^{2}+\beta(t_j) x +
\gamma(t_j) y + 
\delta(t_j)}{(x-\x)(x-\xx)} \,\frac{dx}{y}
\]
where $\alpha\xi_{k}^{2}+\beta \xi_{k} - \gamma 
\eta_{k} + \delta = 0,\;k=1,2$ (here, and henceforth,
we suppress the 
$t_j$). These, together with condition (i) give 
\begin{eqnarray}
\label{eq:alphaetc}
\alpha&=&\frac{\delta}{\x\xx}+\frac{1}{4\pi 
i}\left(\f{\y}{\x}-\f{\yy}{\xx}\right),\\
\beta&=&\frac{\y-\yy}{4\pi i}-(\x+\xx)\alpha,\\
\gamma&=&\frac{\x-\xx}{4\pi i.}
\end{eqnarray}
In terms of $\tilde x$ and $\tilde y$,
\[
{\o}^m=\frac{-(\alpha+\beta \tilde x +\gamma \tilde y + 
\delta \tilde x^2)}{(1-\x \tilde x)(1-\xx \tilde x)}\frac{d\tilde x}{\tilde y}
\]
so from (\ref{eq:tilde_x}), 
\[
\om(P_0)=\f{\eps}{a\cm}(\alpha-\gamma),
\]
and hence 
\[
\dot{\om}(P_0)=\f{\eps}{ac_j}\left[\f{\dot{\delta}}{\x\xx}+
\frac{1}{4\pi i}\left(\dot{\left(\frac{\y}{\x}\right)}-
\dot{\left(\frac{\yy}{\xx}\right)} + \dot{\xx}-\dot{\x}\right)\right].
\]
\begin{lem}
\label{lem:expn_in_dots}
\[
\dot{\delta}_m=0
\] so
\[
\dot{\om}(P_0)=\frac{\eps}{4\pi iac_j}\left(\dot{\left(\frac{\y}{\x}\right)}-
\dot{\left(\frac{\yy}{\xx}\right)} + \dot{\xx}-\dot{\x}\right).
\]
\end{lem}

\begin{proof}
\[
\dot{\o}^m=\frac{\dot{\alpha} x^2 + \dot{\beta} x + \dot{\gamma} y + 
\gamma\dot{y} +\dot{\delta}}{(x-\x)(x-\xx)}\frac{dx}{y}
+ \left(\frac{\dot\x}{x-\x}+\frac{\dot\xx}{x-\xx}-\frac{\dot{y}}{y}\right)\om.
\]
Recalling that when $t_j=0$, $x(P_j)=0$ and $y^2=x^2(x-3a)(x+a)$, we deduce that
\begin{eqnarray*}
0&=&\res_{P_j}\dot{\o}^m\\
&=&\frac{\dot{\delta}}{\sqrt{3}ai}-\res_{P_j}\left(\frac{\dot{y}}{y}\om\right)\\
&=&\frac{\dot{\delta}}{\sqrt{3}ai}+\res_{P_j}\left(\frac{1}{y^2}\om\right).
\end{eqnarray*}
From (\ref{eq:alphaetc}),
\[
\om=\frac{\alpha dx}{y} + \frac{(\y-\yy)x+(\x-\xx)y+\x\yy-\xx\y}{4\pi i (x-\x)(x-\xx)}\frac{dx}{y},
\]
so 
\[
\res_{P_j}\left(\frac{1}{y^2}\om\right)=0
\]
and we conclude that $\dot{\delta}_m=0$.
\end{proof}

\begin{lem}
\label{lem:expn_in_xi_eta}
\[
\dot{\om}(P_0)=\frac{\eps}{12\pi i a c_j}\left( \frac{\yy}{\xx^2}-
\f{\y}{\x^2}+i\sin(3\thmj)\left(\f{1}{\xx^3}-\f{1}{\x^3}\right)+
\f{2a}{\x}-\f{2a}{\xx}\right),
\]
and
\[
\begin{array}{rcl}
\x&=&a(2\cos(\thmj+\2)+1)\\
\xx&=&a(2\cos(\thmj-\2)+1)\\
{\y}/{\x}&=&2ai\sin(\thmj+\2)\\
{\yy}/{\xx}&=&2ai\sin(\thmj-\2)
\end{array}
\]
\end{lem}
\begin{proof}
By definition,
\begin{equation}
x=a\left(\f{\eps\z}{c_j}+1+\f{c_j}{\eps\z}\right), \label{eq:x}
\end{equation}
and using (\ref{eq:y}),
\begin{equation}
\f{y}{x}=a\left(\f{\eps\z}{c_j}-\f{c_j}{\eps\z}\right),\label{eq:y/x}
\end{equation}
so substituting $\z(P_m)=c_m=e^{i\theta_m}$ and $\z(Q_m)=\eps c_m$ gives 
the expressions for $\xi_k$ and $\f{\eta_k}{\xi_k}$, $k=1,2$.

From (\ref{eq:lambda}),
\begin{equation}
\f{d}{dt_j}\,p_{t_j}(\xi_k)=0 \label{eq:deriv:p}
\end{equation}
so 
\[
\dot{\xi}_k=\f{\xi_k-2a}{3\xi_k}.
\]
Differentiating $y$ and using (\ref{eq:deriv:p}) gives
\[
\dot\eta_k=\f{-\eta_k}{3\xi_k}.
\]
Also, we have
\begin{eqnarray*}
y&=&\f{{\lambda}{b_j}^{-1}-{b_j}{\lambda}^{-1}}{2(x-2a)}\\
&=&\f{i\sin{3\thmj}}{x-2a}.
\end{eqnarray*}
Thus
\begin{eqnarray*}
\dot{\left(\f{\eta_k}{\xi_k}\right)}&=&\f{\dot\eta_k}{\xi_k} - 
\f{\dot{\xi_k}\eta_k}{\xi_k^2}\\
&=&\f{-\eta_k}{3\xi_k^2}- \f{i\sin{3\thmj}}{3\xi_k^3}
\end{eqnarray*}
and substituting our expressions into Lemma~\ref{lem:expn_in_dots} 
completes the proof of  Lemma~\ref{lem:expn_in_xi_eta}.  
\end{proof}
\begin{lem} \label{lem:simplification}
\[
\f{\yy}{\xx^2} - \f{\y}{\x^2} = \f{-\sqrt{3}i}{2\cos\thmj + 1}
\]
\end{lem}
The proof is a straightforward but tedious computation, and we omit it.
From Lemma~\ref{lem:expn_in_xi_eta} and Lemma~\ref{lem:simplification} it follows
that
\[
\dot{\om}(P_0) = \f{\eps}{12\pi i a c_j} (\xmj + i \ymj)
\]
where
\[
\xmj=\frac{2a}{\x}-\frac{2a}{\xx}
\]
and
\[
\ymj =
\sin{3\theta_{mj}}\left(\frac{(\x)^3-(\xx)^3}{(\x\xx)^3}\right)
-\frac{\sqrt{3}}{2\cos{\thmj}+1}.
\]
This completes the proof of Proposition~\ref{thm:derivs}.

Since the upper two blocks of $M$ are diagonal, we may make use of the 
following linear algebra lemma, whose proof is left as a straightforward 
exercise for the reader.

\begin{lem}
Suppose $M$ is a  matrix of the form 
$\left(\begin{smallmatrix} A & B\\ C & D\end{smallmatrix}\right)$, 
where $A,\,B,\,C,\,D$ are $n\times n$ 
matrices and $A$, $B$ are diagonal. Then 
$M$ is invertible if and only if $CB-DA$ is invertible.
\end{lem}

\section{Asymptotics}
We have shown that the invertibility of $dW_X$ is equivalent to that of 
the $g\times g$ matrix $N$ with entries
\[
n^m_j= c^m_k b^k_j - d^m_k a^k_j 
\]
where
\[
a^m_j = \left.\frac{\partial u_m}{\partial 
\theta_j}\right|_{\bft=0},\, 
b^m_j = \left.\frac{\partial
v_m}{\partial\theta_j}\right|_{\bft=0},\,
c^m_j = \left.\frac{\partial u_m}{\partial
t_j}\right|_{\bft=0},\,
d^m_j = \left.\frac{\partial v_m}{\partial
t_j}\right|_{\bft=0}.
\]
In particular, notice that $\det(N)$ is an analytic function of
$\theta_1,\ldots,\theta_g$. 

In this section we prove that $dW_X$ is invertible at one of the rational
curves in $\caX$, from which it follows that nearby there are countably many
spectral curves for which $W_X\in Gr(2,g+2)$ is a rational 2-plane i.e.\
there are countably many spectral curves of genus $g$ yielding minimal 
immersed tori in $\CP^2$. A simple modification will prove the analogous 
result for minimal Lagrangian tori in $\CP^{2}$ of spectral  genus $g=2p$.

Our method will be to extend $\det(N)$, thought of as a 
function of $c_1,\ldots,c_g$, to a meromorphic function
$h(z_1,z_2,\ldots,z_g)$ on $(\C^*)^g$ by replacing 
$c_j=e^{i\theta_j}$ by
$z_j=r_je^{i\theta_j}$ in the formulas above. We will then show the following.
\begin{prop} 
\label{prop:asymptotics}
For the complex parameter $z$,
\[
\lim_{z\to\infty}h(z,z^2,\ldots,z^g)\neq 0,
\]
whence $\det(N)$ is not identically zero. Since $\det(N)$ is analytic, it
follows that $dW_X$ is invertible on an open subset of $\caX_b$ about $X_b$. 
For $g=2p$ 
\[
\lim_{z\to\infty} h(z,\ldots,z^p,-z,\ldots,-z^p)\neq 0,
\]
from which it follows that $dW_X$ is invertible on an open subset of
the subspace $\caX_L \subset\caX$ of those Riemann surfaces with an 
involution covering $\lambda\mapsto -\lambda$.
\end{prop}
Theorems \ref{th:every_genus} and \ref{th:every_dim} follow from this.
To begin the proof, we compute explicitly the entries for the matrix $N$.
\begin{lem}
We have
\[
n_j^j = \frac{3}{4\pi^2 a}
\]
and for $m\neq j$
\[
n^m_j = \frac{-\sqrt{3}}{24\pi^2 a} Y_{mj}.
\]
\end{lem}
\noindent
The proof is a straightforward calculation left for the reader.

Now we consider the extension of $Y_{mj}$ to $(\C^*)^g$ as described above.
First we observe, from Lemma~\ref{lem:expn_in_xi_eta} , that
\[
\xi_1 = a[2\cos(\thmj +\frac{2\pi}{3})+1],\quad
\xi_2 = a[2\cos(\thmj -\frac{2\pi}{3})+1].
\]
These are restrictions of the functions
\begin{eqnarray*}
\xi_1 & = & a[\eps z_mz_j^{-1} +\eps^{-1}z_m^{-1}z_j +1]\\
\xi_2 & = & a[\eps^{-1} z_mz_j^{-1} +\eps z_m^{-1}z_j +1]
\end{eqnarray*}
to $\{(z_1,\ldots,z_g):|z_j|=1\}$.  With these expressions the function
$Y_{mj}$ has extension
\begin{equation}
Y_{mj} = \frac{1}{2i}(z_m^3z_j^{-3}-z_j^3z_m^{-3})
\frac{\xi_1^3-\xi_2^3}{\xi_1^3\xi_2^3} - \frac{\sqrt{3}}{z_mz_j^{-1}+z_jz_m^{-1}
+1}.
\end{equation}
Notice that this is a function of $z_mz_j^{-1}$. To prove the proposition we
consider the asymptotics of this expression in two cases: $z_mz_j^{-1}\to
\infty$ and $z_mz_j^{-1}\to 0$. In the first case a computation gives
\[
Y_{mj} = 4\sqrt{3}z_m^{-1}z_j(1+O(z_m^{-1}z_j))
\]
while in the second case we compute
\[
Y_{mj} = -4\sqrt{3}z_mz_j^{-1}(1+O(z_mz_j^{-1})).
\]
Therefore, upon the substitution $z_m=z^m$ we see that since $j\neq m$ we have
\[
\lim_{z\to\infty}Y_{mj} = 0.
\]
It follows that, if $N_0$ denotes the diagonal matrix with diagonal entries
$n_j^j$ then
\[
\lim_{z\to\infty} h(z,\ldots,z^g) = \det(N_0) \neq 0.
\]
This proves part (i) of Proposition~\ref{prop:asymptotics}.
For part (ii) we examine the substitution $z_m=z^m, z_{p+m}=-z^m$ for
$m=1,\ldots,p$. If $m,j\leq p$ or $m,j\geq p$ then the asymptotics are as before.
Otherwise
\[
\lim_{z\to\infty}Y_{mj} = 0,\ \text{for}\ j\neq m\pm p,
\]
while for $j=m\pm p$ we observe that $\xi_1=\xi_2=2a$ so that
\[
Y_{m,m+p}=Y_{j+p,j} = \sqrt{3}.
\]
In this case 
\[
\lim_{z\to\infty} h(z,\ldots,z^p,-z,\ldots,-z^p)=\det(N_1)
\]
where $N_1$ is the block matrix
\[
N_1=\frac{3}{24\pi^2 a}\begin{pmatrix} 6I_p & -I_p\\-I_p & 6I_p \end{pmatrix}
\]
and $I_p$ is the $p\times p$ identity matrix. Since $\det(N_1)\neq 0$ this
completes the proof of Proposition~\ref{prop:asymptotics}.

\end{document}